	\documentstyle[12pt]{article}

\begin{document}

\title{Formulae for the derivatives of   heat semigroups}
\author{K. D.   Elworthy  and X.-M. Li\thanks{Research supported by SERC
grant GR/H67263}}
\date{}
\maketitle

\newcommand{\A}{{\bf \cal A}}
\newcommand{\B}{{ \bf \cal B }}
\newcommand{\C}{{\cal C}}
\newcommand{\F}{{\cal F}}
\newcommand{\G}{{\cal G}}
\newcommand{\h}{{\cal H}}
\newcommand{\K}{{\cal K}}
\newcommand{\half}{{  {1\over 2}  }}
\newcommand{\heatsemif}{{ {\rm e}^{ \half t\triangle^{h,1}}   }}
\newcommand{\heatsemi}{{ {\rm e}^{\half t \triangle^{h}}   }}

\newtheorem{theorem}{Theorem}[section]
\newtheorem{proposition}[theorem]{Proposition}
\newtheorem{lemma}[theorem]{Lemma}
\newtheorem{corollary}[theorem]{Corollary}
\newtheorem{definition}{Definition}[section]

\newtheorem{theorem**}{Theorem}
\newtheorem{theorem*}[theorem**]{Theorem A}
\newtheorem{lemma*}[theorem**]{Lemma A}
\newtheorem{proposition*}[theorem**]{Proposition A}

\def\limsup{\mathop{\overline{\rm lim}}}
\def\liminf{\mathop{\underline{\rm lim}}}
\def\exp{{\rm e}}

\section{Introduction}

Let $M$ be a smooth manifold.  Consider first a non-degenerate stochastic
 differential equation:
\begin{equation} dx_t=X(x_t)\circ dB_t +A(x_t)dt
\label{eq: basic}
\end{equation}
on $M$ with smooth coefficients : $A, X$, where $\{B_t: t\ge 0\}$ is a $R^m$
valued Brownian   motion on a filtered probability space
 $\{\Omega,\F, \F_t, P\}$
. Let $P_t$ be the associated sub-Markovian semigroup and ${\cal A}$
the infinitesimal generator, a second order elliptic operator. In \cite{ELflow} a
formula  for the derivative $d(P_tf)_{x_0}(v_0)$ of $P_tf$ at $x_0$ 
in direction $v_0$ of the form:

\begin{equation}
d(P_tf)_{x_0}(v_0)={1 \over t}{\bf E} f(x_t)\int_0^t < v_s, X(x_s)dB_s>
\label{eq: formula1}
\end{equation}
 was given,  where $v_.$ is a certain stochastic process starting at $v_0$. The
 process $v_.$ could be given either by the derivative flow of (1) or in terms
 of a naturally related curvature. In the latter case and when ${\cal A}
 ={1\over 2}\triangle_M$ for some Riemannian structure the formula reduced to 
one obtained by  Bismut in\cite{Bismut} leading to his well known formula for 
$\nabla \log p_t(x,y)$, the gradient of the logarithm of the fundamental 
solution to the heat equation on a Riemannian manifold. Bismut's proof is in 
terms of Malliavin  calculus while the proofs suggested in \cite{ELflow} 
following the approach of Elliott and Kohlman \cite{ell-koh} are very 
elementary. However there the results were actually given for a compact 
manifold as special cases of a more general result which needed some 
 differential geometric apparatus. Here we show the formula holds in  a more 
general context, extend it  to higher derivatives,  and give  similar formulae
 for differential forms of all orders  extracted from \cite{Li.thesis}.  In 
particular we have a  simple  proof of the formulae for somewhat more  general
 stochastic differential equations.

One  importance of these formulae is that they demonstrate the smoothing
effect of $P_t$ showing clearly what happens at $t=0$. To bring out
the simplicity we first give proofs of the basic results for It\^o equations
on $R^n$.

 There are extensions to infinite  dimensional systems with applications to 
smoothing and the strong Feller property  for infinite dimensional Kolmogorov
 equations in \cite{DA-EL-ZA} \cite{PE-ZA}.
There are also applications to non-linear reaction-diffusion equations\cite
{Li-Zhao}.  For other generalizations of Bismut's formula in a geometric
 context see  \cite{Norris93}. The work of Krylov \cite{Krylov93} in
this general area  must also be mentioned although the approach
and aims  are rather different.

Throughout this article, we shall use $BC^r$ for the space of bounded 
$C^r$ functions  with their first $r$ derivatives bounded (using a given
 Riemannian metric on  the manifold).

\section{Formulae with simple proof for $R^n$}

For $M=R^n$, we can take the It\^o form of (1):
\begin{equation} dx_t=X(x_t) dB_t +Z(x_t)dt
\label{eq: Ito}
\end{equation}
where $X\colon R^n \to L(R^m, R^n)$ and $Z\colon R^n \to R^n$ are $C^\infty$ 
with derivatives $DX\colon$  $R^n$ $\to L(R^n, L(R^m,R^n))$ and 
$DZ\colon R^n \to L(R^n, R^n)$ etc. There is the derivative equation:

\begin{equation}
dv_t=D X(x_t)(v_t)dB_t +DZ(x_t)(v_t)dt.
\label{eq: covIto}
\end{equation}

\noindent
whose solution $v_t=DF_t(x_0)(v_0)$ starting from  $v_0$ is the derivative 
(in probability) of $F_t$ at $x_0$ in the direction $v_0$. Here 
$\{F_t(-), t\ge 0\}$ is a solution flow to (3), so that $x_t=F_t(x_0)$, for 
$x_0\in R^n$. We do not need to assume the existence of a sample smooth version
 of $F_t\colon M\times \Omega \to M$.

\bigskip

For $\phi\colon R^n \to L(R^n; R)$, define $\delta P_t(\phi)\colon R^n \to
 L(R^n, R)$ by

\begin{equation}
\left (\delta P_t(\phi)\right)_{x_0}(v_0)=E\phi_{x_t}\left(v_t\right)
\end{equation}

\noindent
whenever the right hand side exists. Here $\phi_x(v)=\phi(x)(v)$. In particular this can be applied to $\phi_x=(df)_x\colon= Df(x)$ where $f\colon R^n \to R$ has bounded  derivative. Formal differentiation under the expectation suggests 
$$d(P_tf)_{x_0}(v_0)=(\delta P_t(df))_{x_0}(v_0).$$
 This is well known when $X$ and $Z$ have bounded first derivatives.
 It cannot hold for $f\equiv 1$
  when (3) is not complete(i.e. explosive). In fact we will
 deal  only with  complete systems: we are almost forced to do this since for
 $\delta P_t$ to have a reasonable domain of definition some integrability 
conditions on $DF_t(x_0)$ are needed and it is shown in \cite{Li.thesis} that
 non-explosion follows,  for a  wide class of symmetrizable diffusions, from
 $dP_tf=\delta P_t(df)$ for all $f\in C_K^\infty$ together with 
$E\chi_{t<\xi}|DF_t(x_0)|<\infty$ for all $x_0\in M, t>0$. Here $\xi$ is the
explosion time. Precise conditions 
for $d(P_tf)=(\delta P_t)(df)$ are given  in an appendix  below.

\bigskip

\noindent
Our basic result is the following. It originally appeared in this form in
\cite{Li.thesis}.

\begin{theorem}\label{th: Elworthy-Li} Let (3) be complete and non-degenerate,
so there is a right inverse map $Y(x)$ to   $X(x)$  for each $x$ in $R^n$,
 smooth in $x$.  Let $f\colon R^n \to R$ be $BC^1$ with  $\delta P_t(df)=
d(P_tf)$ almost surely (w.r.t. Lebesgue measure) for $t\ge 0$. Then for 
 almost all $x_0\in R^n$ and $t>0$

\begin{equation}
d(P_tf)(x_0)(v_0)={1\over t} Ef(x_t) \int_0^t <Y(x_s)(v_s),dB_s>_{R^m},
\hskip6pt v_0\in R^n
\label{EL-LI}
\end{equation}

\noindent
provided $\int_0^t <Y(x_s)(v_s),dB_s>_{R^m}, t\ge 0$ is a martingale.
\end{theorem}

\noindent {\bf Proof:}
Let $T>0$.  Parabolic regularity ensures that  It\^o's formula can be applied
 to  $(t,x) \mapsto P_{T-t}f(x), 0\le t\le T$ to yield:
 
\begin{equation}
P_{T-t}f(x_t)=P_Tf(x_0)+ \int_0^t d (P_{T-s}f)_{x_s}(X(x_s)dB_s).
\end{equation}

\noindent
for $t\in [0,T)$. Taking the limit as $t\to  T$, we have:
$$f(x_T)=P_Tf(x_0)+\int_0^T d (P_{T-s}f)_{x_s} (X(x_s)dB_s).$$
Multiplying through by our martingale and then taking expectations using the 
fact that $f$ is bounded, we obtain:

\begin{eqnarray*}
&&Ef(x_T)\int_0^T<Y(x_s)v_s, dB_s>=
E\int_0^Td (P_{T-s}f)_{x_s}(v_s)ds\\
&=& E\int_0^T\left((\delta P_{T-s})(df)\right)_{x_s}(v_s)ds
= \int_0^T\left((\delta P_s)\left((\delta P_{T-s})(df)\right)\right)_{x_0}(v_0)ds \\
&=&\int_0^T (\delta P_T(df))_{x_0}(v_0)ds
=T\ \delta P_T(df)_{x_0}(v_0).
\end{eqnarray*}
by the equivalence of the law of $x_s$ with Lebesque measure and the semigroup
property of $\delta P_t$.
\hfill \rule{3mm}{3mm}

\bigskip

\noindent
{\underline {\bf Remarks:}}

1. The proof shows  that under our conditions  equality in (\ref{EL-LI}) holds
 for each $x_0\in M$ if and only if $\delta P_t(df)=d(P_tf)$ at each point.
 This is true provided  $x\mapsto E|DF_t(x)|$  is continuous. The same will
hold for various variations of theorem \ref{th: Elworthy-Li} which follow.

2.  The martingale hypothesis is satisfied if 
$$\int_0^t E|Y(x_s)(v_s)|^2ds<\infty$$ 
for all $t$. In turn this is implied
by the uniform ellipticity condition

\noindent  $|Y(x)(w)|^2\le {1\over \delta } |w|^2$
for all $x, w\in R^n$, for some $\delta>0$, together with
\begin{equation}
\int_0^t E|v_s|^2ds<\infty, \hskip 12pt t\ge 0.
\end{equation}

\noindent
Under these conditions, (\ref{EL-LI}) yields
$$\sup_{x\in R^n} |d(P_tf)_x|\le {1\over t} \sup_{x\in R^n} |f(x)|{1\over \delta}
\sup_{x\in R^n}\sqrt{\int_0^tE|DF_s(x)|^2}ds.$$

\noindent
In particular if $X$, $Z$ have bounded first derivatives, then Gronwall's
 inequality together with
(\ref{eq: covIto}) yields a constant  $\alpha$ with
\begin{equation}
\sup_{x\in R^n} |d(P_tf)_x|  \le {1\over \delta} {1\over \alpha t}
\sqrt{\exp^{\alpha t}-1}\sup_{x\in R^n}|f(x)|.
\label{eq:1}
\end{equation}

For Sobolev norm estimates see (\ref{eq: Sobolev 100}) below.

\begin{corollary}
Let (1) be complete and uniformly elliptic.  Then (\ref{EL-LI})  holds for
all $f$ in $BC^1$ provided that $H_2(x)(v,v)$ is bounded above, i.e.
$H_2(x)(v,v)\le c|v|^2$. Here $H_2$ is defined by:

$$H_2(x)(v,v)=2<DZ(x)(v),v>+\sum_1^m|DX^i(x)(v)|^2+
\sum_1^m{1\over |v|^2} <DX^i(x)(v),v>^2$$

\end{corollary}

\noindent {\bf Proof:}
By lemma A\ref{le: 7.2}, we have $\int_0^tE|v_s|^2ds$ finite
for each $t>0$ while theorem A\ref{th: 7.5}  and its remark
gives us the a.s. differentiability required.
\hfill\rule{3mm}{3mm}

\bigskip

\noindent
2. The case when there is a zero order term and when the coefficients
are time dependent can be dealt with in the same way: Let
$\{\A_t:t\ge 0\}$ be second order elliptic operators on $R^n$ with
\begin{eqnarray*}
\A_t(f)(x)&=&\half {\rm trace}D^2f(x)\left(X_t(x)(-),X_t(x)(-)\right)\\
&+&Df(x)\left(Z_t(x)\right)+V_t(x)f(x)
\end{eqnarray*}
for $X_t, Z_t$ as $X, Z$ before, for each $t>0$  continuous in $t$
 together with their  spatial derivatives, and with $V_\cdot(\cdot)
\colon [0,\infty)\times R^n\to R$ continuous and bounded above on each
 $[0,T]\times R^n$.
For each $T>0$ and $x_0\in R^n$ let $\{x_t^T\colon 0\le t\le T\}$ be the
solution of

$$dx_t^T=X_{T-t}(x_t^T)dB_t+Z_{T-t}(x_t^T)dt,$$
with $x_0^T=x_0$ (assuming no explosion) and set
$$\alpha_t^T(x_0)=\exp^{\int_0^tV_{T-s}(x_s^T)ds} \hskip 24pt
0\le t\le T.$$
Also write $x_t^T(\omega)=F_t^T(x_0,\omega)$.
Now suppose $u_\cdot(\cdot): [0,\infty)\times R^n\to R$ satisfies:
\begin{equation}   {\partial u_t \over \partial t}=\A_tu_t,
\hskip 24pt t>0 \end{equation}
and is $C^{1,2}$ and  bounded on each $[0,T]\times R^n$. Then, as before,
 we can apply
It\^o's formula to $\{u_{T-t}(x): 0\le t\le T\}$ to see that 
$\{u_{T-t}(x_t^T): 0\le t\le T\}$ is a martingale and 
$ u_t(x)=E\alpha_t^t(x_0)u_0(x_t^t)$, e.g. see \cite{Freidlin85}.
If we also assume:

\noindent (i) $V_t$ is $C^1$ for each $t$ and continuous and bounded above on
each $[0,T]\times R^n$.

\noindent
 (ii)  We can differentiate under the expectation to have,  for almost
all $x_0\in R^n$
$$Du_t(x_0)(v_0)=E\left( \alpha_t^t(x_0)Du_0(x_t^t)(v_t^t)
 + \alpha_t^t(x_0)u_0(x_t^t)\int_0^tDV_{t-s}(x_s^t)v_s^tds\right)$$

\noindent
where $v_s^t$ solves
\begin{eqnarray*}
dv_s^t&=&DX_{t-s}(x_s^t)(v_s^t)dB_s +DZ_{t-s}(x_s^t)(v_s^t)ds\\
v_0^t&=&v_0 \hskip 24pt 0\le s\le t.
\end{eqnarray*}

\noindent
(iii) For $Y_t(x)$ a right inverse for $X_t(x)$, assume
$\int_0^t<Y_{T-s}(x_s^T)(v_s^T), dB_s>$,  $0\le t\le T$, is a martingale.

 Then  for each $0<t\le T$:

\begin{equation}\label{eq: Remark 10} \begin{array}{ll}
Du_t(x_0)(v_0)&={1\over t}Eu_0(x_t^t)\exp^{\int_0^t V_{t-s}(x_s^T)ds}
\int_0^t<Y_{t-s}(x_s^t)(v_s^t),dB_s>\\
&+{1\over t}E u_0(x_t^t)\exp^{\int_0^t V_{t-s}(x_s^T)ds}\int_0^t
\int_0^s DV_{t-r}(x_r^t)(v_r^t)drds.
\end{array}\end{equation}

\noindent
The only real additional ingredients in the proof are the almost sure
identities:
$$F_r^{T-s}\left(F_s^T(x_0,\omega), \theta_s(\omega)\right)
=F_{s+r}^T(x_0,\omega),
\hskip 20 pt (x_0,\omega)\in M\time \Omega$$
and
$$\alpha_s^T(x_0,\omega)\alpha_{T-s}^{T-s}\left(F_s^T(x_0,\omega), \theta_s(\omega)\right)
=\alpha_T^T(x_0,\omega).$$
where $\theta_s\colon \Omega \to \Omega$ is the shift, e.g. using the
canonical representation of $\{B_t:t\ge 0\}$.

\bigskip

Note that for $X, Z$ with first two derivatives bounded  and $f$
 in $BC^2$,  we can differentiate twice
under the integral sign \cite{ELbook} to see directly that $P_{T-t}f(x)$ is
sufficiently regular to prove (\ref{EL-LI}).  This gives (\ref{EL-LI})
 without using elliptic regularity results and from this (e.g. via  (\ref
{eq:1})) we can  approximate to obtain the smoothing property directly (see 
\cite{DA-EL-ZA} for this approach in infinite dimensions). For further
 smoothing, we can use the next result:  ($c$ is a constant)

\begin{theorem}\label{th: 2nd derivative}
 Assume that equation (3) is complete and has 
uniform ellipticity: $X$ has a right inverse $Y$, which is bounded on $R^n$.
Suppose also
\begin{enumerate}
\item  For each  $x_0, u_0\in R^n$ and  each $T>0$:
\begin{equation}
\int_0^T E|DF_s(x_0)(u_0)|^2 ds\le c |u_0|^2,
\label{eq: hypothesis 1}
\end{equation}

\item For each $t>0$,

$$\sup_{0\le s\le t}\sup_{y_0\in R^n} 
 \left(E|D^2F_s(y_0)(u_0,v_0)| \right)\le c|u_0||v_0|,$$ 
and
$$\sup_{0\le s\le t}\sup_{y_0\in R^n}\left(E|DF_s(y_0)|\right)
\le c.$$
\end{enumerate}
Let  $f$ be in $BC^2$  and such that $d(P_tf)_{x_0}=\delta P_t(df)_{x_0}$  for
almost all $x_0\in R^n$   and that   we can differentiate $P_tf$ under the
 expectation to give, for almost all $x_0$:

\begin{equation}
\begin{array}{lll}
D^2P_tf(x_0)(u_0)(v_0)&=& ED^2f(x_t)(DF_t(x_0)u_0, DF_t(x_0)v_0)\\
&+& EDf(x_t)(D^2 F_t(x_0)(u_0, v_0))\\
\end{array} 
\label{eq: 1a} \end{equation}

\noindent
for each $t\ge 0$. Then for almost all $x_0$ in $R^n$ and all $t>0$,

\begin{equation}
\begin{array}{ll}
D^2P_tf(x_0)(u_0,v_0)&=
{4\over t^2} E\left\{f(x_t)\int_{t\over 2}^t   <Y(x_s)v_s, dB_s> \int_0^
{t\over 2}<Y(x_s)u_s, dB_s>\right\}\\
&-{2\over t}E\int_0^{t\over 2} D(P_{t-s}f)(x_s)(DX(x_s)(v_s)(Y(x_s)u_s))ds\\
&+{2\over t}E\int_0^{t\over 2}   D(P_{t-s}f)(x_s)(D^2F_s(x_0)(u_0,v_0))ds.\\
\end{array}
\label{eq:2a}
\end{equation}

\noindent 
 If also  $\int_0^{t\over 2}<DY(x_s)(DF(x_0)u_0)(DF(x_0)v_0), dB_s>$
is a martingale,
 then

\begin{equation}
\begin{array}{ll}
 & D^2P_tf(x_0)(u_0,v_0)\\
&={4\over t^2}E\left\{f(x_t)\int_{t\over 2}^t <Y(x_s)v_s, dB_s>
\int_0^{t\over 2}<Y(x_s)u_s, dB_s>\right\}\\
&+{2\over t}E\left\{f(x_{t})\int_0^{t\over 2}
<DY(x_s)(u_s)(v_s), dB_s>\right\}\\
&+ {2\over t}E\left\{ f(x_t)\int_0^{t\over 2} 
<Y(x_s)D^2F_s(x_0)(u_0,v_0), dB_s>\right\}.\\
\end{array}
\label{eq:2}
\end{equation}
\end{theorem}

\noindent{\bf Proof:}
Since  $d(P_{T-t}f)$  is smooth and satisfies 
the relevant parabolic equation,  by It\^o's formula (e.g.\cite{ELflow} 
cor. 3E1), if $0\le t<T$,

\begin{eqnarray*}
d(P_{T-t} df)_{x_t}(v_t)&=&d (P_T f)_{x_0}(v_0)+
\int_0^t\nabla\left(d( P_{T-s}f)_{x_s}\right)(X(x_s)dB_s)(v_s)\\
&+&\int_0^t(d( P_{T-s}f)_{x_s}\left(DX(x_s)(v_s)dB_s\right )\\
\end{eqnarray*}

\noindent  giving

\begin{eqnarray*}
(df)_{x_T}(v_T)&=& d( P_Tf)_{x_0}(v_0)+\int _0^T D^2(P_{T-s}f)(x_s)
(X(x_s)dB_s)(v_s)\\
&+&\int_0^T   D(P_{T-s}f)(x_s)(DX(x_s)(v_s)dB_s).\\
\end{eqnarray*}
Using  the uniform ellipticity and hypothesis 1 (i.e. equation 
(\ref{eq: hypothesis 1})) this gives

\begin{eqnarray*}
&&E\left\{(df)_{x_T}(v_T)\int_0^T <Y(x_s)u_s, dB_s>\right\}\\
&&=E\int_0^T D^2(P_{T-s}f)(x_s)(u_s)(v_s)ds\\
&&+E\int_0^T  D(P_{T-s}f)(x_s)(DX(x_s)(v_s)(Y(x_s)u_s))ds.\\
\end{eqnarray*}

\noindent
Thus by ($\ref{eq: 1a}$), and using the two hypotheses to justify  changing 
the order of integration,
\begin{eqnarray*}
T\left[ D^2P_Tf(x_0)(u_0,v_0)\right]
&&=E\left\{Df(x_T)(v_T)\int_0^T<Y(x_s)u_s, dB_s>\right\}\\
  &-&E\left\{\int_0^TD(P_{T-s}f)(x_s)(DX(x_s)(v_s)(Y(x_s)u_s))ds\right\}\\
  &+& E\int_0^T\left\{D(P_{T-s}f)(x_s)
\left(D^2F_s(x_0)(u_0,v_0)\right)\right\}ds.\\
\end{eqnarray*}

\noindent
Now let $T=t/2$, and replace $f$ by $P_{t\over 2}f$. Note that by theorem 
~\ref{th: Elworthy-Li}  and the Markov property(or cocycle property of flows)

$$DP_{t\over 2}f(x_{t\over 2})(v_{t\over 2})=
{2\over t} E\left\{f(x_t)\int_{t\over 2}^t
<Y(x_s)v_s, dB_s>| x_s: 0\le s\le t/2\right\}.$$

\noindent  We see
\begin{eqnarray*}
& &D^2 P_tf(x_0)(u_0,v_0)\\
&=&{4\over t^2}E\left\{ f(x_t)\int_{t\over 2}^t <Y(x_s)v_s, dB_s>
\int_0^{t\over 2}<Y(x_s)u_s, dB_s>\right\}\\
&-&{2\over t} E\int _0^{t\over 2}D(P_{t/2-s}f)(x_s)(DX(x_s)(v_s)(Y(x_s)u_s))ds\\
&+&{2\over t} E\int_0^{t\over 2} D(P_{t/2-s}f)(x_s)(D^2F_s(x_0)(u_0,v_0))ds,\\
\end{eqnarray*}

\noindent giving (\ref{eq:2a}).
Now apply It\^o's formula to $\{P_{t-s}f(x_s): 0\le s <t\}$ at $s={t\over 2}$
to obtain
\begin{equation}
P_{t\over 2}f(x_{t\over 2})=P_tf(x_0)+
\int_0^{t\over 2} D(P_{t/2-s}f)(x_s)(X(x_s)dB_s).
\label{eq:8}
\end{equation}

The equation (\ref{eq:2}) follows on multiplying (\ref{eq:8}) by
$\int_0^{t\over 2} <DY(x_s)(u_s)(v_s), dB_s>$  and also by 
$\int_0^{t\over 2} <Y(x_s)D^2F_s(x_0)(u_0,v_0), dB_s>$ and taking expectations
 to replace the 2nd and 3rd terms in the right hand side of (\ref{eq:8}),
using the identity:

\begin{equation}
DX(x)(u)(Y(x)v)+X(x)DY(x)(u)(v)=0.
\label{eq:9}
\end{equation}

\noindent{\bf \underline{Remarks} }

(A).  Formulae (\ref{eq:2a}) combined with theorem \ref{th: Elworthy-Li} has 
some  advantage over (\ref{eq:2}) for estimation since the derivative of $Y$
 does not appear.

(B). Formulae (\ref{eq:2}) can be obtained by applying theorem 
~\ref{th: Elworthy-Li}, with $t$ replaced by $t/2$, to $P_{t\over 2}f$ and then
 differentiating under the expectation and stochastic integral sign, assuming
 this  is legitimate, then using the Markov property to replace the 
$P_{t\over 2}f(x_{t\over 2})$ by $f(x_t)$.

(C). The hypotheses 1 and 2 of the theorem and the conditions on the function
 $f$ are satisfied if  $|DX|$, $|D^2X|$,
  $ |DA|$, and $|D^2 A|$ are bounded.  See lemma A\ref{le: 7.2},
theorem A\ref{th: 7.5} and proposition A\ref{pr: 7.8}.
Furthermore the martingale condition needed for (\ref{eq:2}) also holds if
 $DY$ is bounded as a bilinear map.

\section{Formulae with simple proof for $M$}

For a general smooth manifold $M$, we return to the Stratonovich equation (1). 
We will continue to assume non-explosion and non-degeneracy.  Thus now $X(x)$ 
is a surjective linear map of $R^m$ onto the tangent space $T_xM$ to $M$ at $x$
 and $A$ is a  smooth vector field on $M$. Write $X^i(x)=X(x)(e_i)$ for
 $e_1,\dots. e_m$ an orthonormal basis for $R^m$. Thus (1) becomes:

\begin{equation}
dx_t=\sum_1^m X^i(x_t)\circ dB_t^i +A(x_t)dt
\label{eq: basic1a}
\end{equation}

\noindent
Here $\{B_t^i, t\ge 0\}$ are independent one dimensional Brownian motions.
The generator $\A$, being elliptic, can be written
 $\A=\half \triangle +Z$ where
$\triangle$ denotes the Laplace Beltrami operator for an induced Riemannian
metric on $M$ and $Z$ is a smooth vector field on $M$. Using this metric and
the Levi-Civita connection

\begin{equation}
Z=A^X=\half \sum_1^m \nabla X^i(X^i(x))+A
\label{eq:sec3 2}
\end{equation}

\noindent
The derivative equation extending (\ref{eq: covIto}) is most concisely expressed as a covariant equation

\begin{equation}
dv_t=\nabla X(v_t)\circ dB_t +\nabla A(v_t)dt.
\label{eq:cov}
\end{equation}

By definition, this means:
\begin{equation}
d\tilde{v_t}=//_t^{-1}\nabla X\left(//_t\tilde{v_t}\right)\circ dB_t +//_t^{-1}
\nabla A(//_t \tilde v_t)dt
\end{equation}

\noindent
for $\tilde{v_t}=//_t^{-1}v_t$ with $//_t\colon T_{x_0}M\to T_{x_t}M$ parallel
translation along the paths of $\{x_t\colon t\ge 0\}$.

 Recall that  covariant differentiation gives linear maps:

$$\nabla A\colon T_xM\to T_xM, \hskip 20pt x\in M,$$
$$\nabla X\colon T_xM \to L(R^m; T_xM) \hskip 15pt x\in M$$
and
$$\nabla^2 A\colon  T_xM \to L(T_xM; T_xM) \hskip 13pt x\in M$$
sometimes considered as a bilinear map by
$$\nabla^2 A(u,v)=\nabla ^2 A(u)(v) \hskip 6pt etc.$$

For the (measurable) stochastic flow $\{F_t(x)\colon t\ge 0, x\in M\}$ to (1), 
the derivative in probability now becomes a linear map between tangent spaces
 written
$$T_{x_0}F_t\colon T_{x_0}M\to T_{x_t}M  \hskip 20pt x_0\in M,$$
or
$$TF_t\colon TM \to TM,$$
and $v_t=T_{x_0}F_t(v_0)$, the derivative at $x_0$ in the direction $v_0$.

\bigskip

Analogous to the probability semigroup $P_t$, there is the following semigroup 
(formally) on differential forms:

\begin{equation}\label{eq: definition of semigroups}
\delta P_t\phi(v_1,\dots, v_p)= E\phi(TF_t(v_1), \dots, TF_t(v_p)).
\end{equation}

\noindent
Here $\phi$ is a $p$-form.  If $\phi=df$ for some function $f$, then
$$\delta P_t(df)(v)=Edf(TF_t(v)).$$

In \cite{ELbook}, it was shown  that $\delta P_t(df)=d(P_tf)$  if $\nabla X$, 
$\nabla A$, and $\nabla^2 X$ are bounded, and  if  the stochastic differential
equation is strongly complete on $R^n$( or on  a complete Riemannian manifold
  with bounded curvature).  Theorems of this kind  are since much improved
 partially due to the  concept of strong 1-completeness \cite{Li.thesis}.
 See the appendix  for the definition of strong 1-completeness.

\bigskip

To differentiate $P_tf$ twice it is convenient to use the covariant derivative
$\nabla TF_t$ which is bilinear
$$\nabla T_{x_0}F_t\colon T_{x_0}M \times T_{x_0}M\to T_{x_t}M.$$
It can be defined by

\begin{equation}
\nabla T_{x_0}F_t(u_0,v_0)={D\over \partial s} T_{\sigma(s)}F_t(v(s))|_{s=0}
\end{equation}
for $\sigma$ a $C^1$ curve in $M$ with $\sigma(0)=x_0, \dot{\sigma}(0)=u_0$
 and for $v(s)$ the parallel translate of $v_0$ along $\sigma$ to $\sigma(s)$,
the derivative being a derivative in probability in general, \cite{ELbook} 
  page 141.
 
\bigskip

\noindent
The extensions of theorems \ref{th: Elworthy-Li} and \ref{th: 2nd derivative}
 can be written as follows and proved in essentially the same way; note that
 we can take $Y(x)=X(x)^*$:

\begin{theorem} Let $M$ be a complete Riemannian manifold and 
$\A=\half \triangle +Z$.
Assume (1) is complete. Let $f\colon M\to R$ be $BC^1$ with:
 $\delta P_t(df)=d(P_t f)$ a.e. for $t\ge 0$. Then for almost all $x_0\in M$,

\begin{equation}\label{eq: formula on manifold}
dP_tf(v_0)={1\over t}Ef(x_t)\int_0^t <v_s, X(x_s)dB_s>_{x_s}, \hskip 6pt 
 v_0\in T_{x_0}M
\end{equation}
provided $\int_0^t <v_s, X(x_s)dB_s>$ is a martingale.
Furthermore assume: 
\begin{enumerate}
\item
  For each $T>0$ and $x_0\in M$,
\begin{equation}\label{eq: manifold hypothesis 1}
\int_0^T E|T_{x_0}F_s(u_0)|^2 ds\le c |u_0|^2, \hskip4pt u_0\in T_{x_0}M,
\end{equation}

\item  For each $T>0$,
\begin{equation} \label{eq: manifold hypothesis 2}
\sup_{0\le s\le T}\sup_{y_0\in M} 
 \left(E|\nabla T_{y_0}F_s(u_0,v_0)| \right)\le c|u_0||v_0|,
\end{equation}
and
\begin{equation}\label{eq: manifold hypothesis 3}
\sup_{0\le s\le T}\sup_{y_0\in M}\left(E|T_{y_0}F_s|\right)
\le c.
\end{equation}
\end{enumerate}

Let $f$ be a $BC^2$ function such  that  we can differentiate $P_tf$ under the expectation to give:

\begin{equation}
\begin{array}{ll}
\nabla dP_tf(-)(-)=& E\nabla df\left(TF_t(-), TF_t(-)\right)\\
&+ Edf(\nabla TF_t(-,-)), \hskip 6pt a.e.\\
\end{array} 
\label{eq:sec3 4} \end{equation}

 \noindent
for each $t\ge 0$. Then for  almost all $x_0\in M$, all   $u_0$, $v_0$ in
 $T_{x_0}M$, and  $t>0$
\begin{eqnarray*}
&&\nabla d(P_tf)(u_0,v_0)\\
&=&{4\over t^2}E\left\{f(x_t)\int_{t\over 2}^t <v_s, X(x_s)dB_s>
\int_0^{t\over 2}<u_s, X(x_s)dB_s>\right\}\\
&+&{2\over t}E\left\{f(x_{t})\left(\int_0^{t\over 2}
<v_s, \nabla X(u_s)dB_s>+
\int_0^{t\over 2} <\nabla TF_s(u_0,v_0), X(x_s)dB_s>\right)\right\}.\\
\label{eq:sec3 5}
\end{eqnarray*}
\label{th: 3.1}
\hfill\rule{3mm}{3mm}
\goodbreak
\end{theorem}

From  theorem  \ref{th: 3.1}  formula $(\ref{eq: formula on manifold})$
  holds for all $x$ if $H_2(x)(v,v)\le c|v|^2$
for some constant $c$,  by lemma A\ref{le: 7.2}, theorem A\ref{th: 7.5}
and its remark. Here

\begin{equation}\label{eq: H2}\begin{array}{ll}
H_2(x)(v,v):=&-Ric_x(v,v)+2<\nabla Z(x)(v),v>+\sum_1^m |\nabla X^i(x)|^2\\
&+\sum_1^m {1\over |v|^2}<\nabla X^i(x)(v), v>^2.
\end{array}
\end{equation}
Suppose the first three derivatives of  $X$ and the first two of $A$ are
bounded,  then all the conditions of the theorem hold.
 See lemma A\ref{le: 7.2}, proposition A\ref{pr: 7.6},  and proposition
A\ref{pr: 7.8}  for details.

\bigskip

Now let $p_t\colon M\times M \to R$, $t>0$ be the heat kernel, (with respect
 to  the Riemannian volume element) so that

\begin{equation}
P_tf(x)=\int_M p_t(x,y)f(y)dy.
\label{eq: heat kernel}
\end{equation}

\bigskip

There is the following Bismut type formula (see \cite{ELflow} and section 5A
 below).

\begin{corollary}\label{co:  3.2}
Suppose  $\delta P_t(df)=d(P_tf)$   for all $f$ in $C_K^\infty$ and for
 all  $t> 0$. Then, for $t>0$,

\begin{equation}
\nabla \log p_t(\cdot, y)(x_0)
={1\over t} E\{\int_0^t \left(TF_s\right)^*X(x_s)dB_s |x_t=y\}
\label{eq: Bismut}
\end{equation}
for almost  all $y\in M$ provided $\int_0^t <v_s, X(x_s)dB_s>$ is a martingale.
In particular  (\ref{eq: Bismut}) holds if $H_2$ defined in (\ref{eq: H2})
is bounded above.
\end{corollary}

\noindent
{\bf Proof:}  
The proof is just as for the compact case. Let $f\in C_K^\infty$. 
 By the smoothness of $p_t(-,-)$ for $t>0$, we can differentiate equation 
 (\ref{eq: heat kernel}) to obtain:

\begin{equation} d(P_tf)(v_0)=\int_M<\nabla p_t(-,y),v_0>_{x_0}f(y)dy.
\label{eq: manifold 100}\end{equation}

\noindent
On the other hand, we may rewrite  equation $(\ref{eq: formula on manifold})$
 as follows:

$$d(P_tf)(v_0)=\int_M  p_t(x_0,y)f(y)E\left\{{1\over t} \int_0^t <TF_s(v_0), X(x_s)dB_s>| x_t=y\right\} dy$$
Comparing the last two equations, we get:
$$\nabla p_t(-,y)(x_0)=p_t(x_0,y)E\left\{{1\over t} \int_0^t TF_s^*(XdB_s)| x_t=y\right\}.$$ 
 \hfill \rule{3mm}{3mm}\goodbreak

\bigskip

Equality in (\ref{eq: Bismut}) for all $y$ will follow from the continuity of
the right hand side in $y$: for this see \cite{Bismut}, the Appendix to
\cite{Norris93}, or \cite{Watanabe}.

\bigskip

Let $h\colon M \to R$ be a smooth function. There is the corresponding 
Sobolev space $W^{p,1}=\{f\colon M\to R$ s.t.
 $f, \nabla f\in L^p(M, \exp^{2h} dx)\}$ for $1\le p\le \infty$ with norm
 $|f|_{L^{p,1}}=|f|_{L^p}+|\nabla f|_{L^p}$
. Here $dx$ is the Riemannian volume measure. 

\begin{corollary} Suppose $\A=\half \triangle+\nabla h$ for smooth $h$   and
that   $$k^2=:\sup_{x\in M}   E  \int_0^t  |T_xF_s|^2  ds <\infty.$$
Then (\ref{eq: formula on manifold}) holds almost everywhere for any
 $f\in L^p$, $1< p\le \infty$, and for $t>0$, $P_t$  gives a continuous map
$$P_t\colon L^p(M, \exp^{2h}dx)\to W^{p,1}(M, \exp^{2h}dx),
 \hskip 12pt 1<p\le \infty$$
with
 \begin{equation}\label{eq: Sobolev 100}
 |(P_tf)|_{L^{p,1}}\le  (1+{k_p\over t}) |f|_{L^p}, \end{equation}
where $k_p=k$ for $2\le p\le \infty$, and $k_p=c_pk^p$ for $1<p<2$ and
$c_p$ a universal constant.
\end{corollary}

\noindent {\bf Proof:}
Take $f$ in $ BC^1$. Noting that $\exp^{2h}dx$ is an invariant measure for the
 solution  of (1),  formula (\ref{eq: formula on manifold}) gives:
\begin{eqnarray*}
|\nabla (P_tf)(v)|_{L^2}&\le& {1\over t}
\sqrt{\int_M \left[Ef(F_t(x))\int_0^t <X(F_s(x))dB_s, TF_s(v)>\right]^2
\exp^{2h} dx}\\
&\le&{1\over t}\left(\sup_{x\in M} E\int_0^t|T_xF_s(v)|^2ds\right)^\half
\sqrt{\int_M Ef(F_t(x))^2\exp^{2h}  dx}\\
&=& {1\over t}\left(\sup_{x\in M} E\int_0^t|T_xF_s(v)|^2ds\right)^\half
|f|_{L^2}.
\end{eqnarray*}

If $f\in L^2$, let $f_n$ be a sequence in $C_K^\infty$ converging to $f$ in
 $L^2$, then $d(P_tf_n)$ 
converges in $L^2$ by the above estimate with limit $ d(P_tf)$.
So formula (\ref{eq: formula on manifold}) holds almost  everywhere for 
$L^2$ functions.

\noindent
On the other hand  if $f$ also belongs to $L^\infty$,
\begin{equation}\label{eq: Sobolev 200}
|P_tf|_{L^{\infty, 1}} \le \left( 1+ \sup_{x\in M} {1\over t}
\left(\int_0^t E|T_xF_s|^2 ds\right)^\half\right) |f|_{L^\infty}.
\end{equation}
 By the Reisz-Thorin interpolation theorem, we see for  $f\in L^2\cap L^p$,

\noindent
 $2\le p\le \infty$,

\begin{equation}
 |(P_tf)|_{L^{p,1}}\le \left(1+{k\over t}\right)|f|_{L^p}.
\end{equation}

 Again we conclude that 
(\ref{eq: formula on manifold}) holds for $f\in L^p$, $2\le p< \infty$.
 For

\noindent  $1<p<2$, let $q$ be  such that 
 ${1\over p}+{1\over q}=1$. Then   H\"older's inequality gives:
\begin{eqnarray*}
&&|\nabla (P_tf)(v)|_{L^p}\le {1\over t}
\left(\int_M \left[Ef(F_t(x))\int_0^t <X(F_s(x))dB_s, TF_s(v)>\right]^p
\exp^{2h} dx \right)^{1\over p}\\
&\le&{1\over t}\left(\sup_{x\in M}
 E\left[\int_0^t <XdB_s, T_xF_s(v)>\right]^q\right)^{1\over q}
\left(\int_M E\left[f(F_t(x))\right]^p\exp^{2h} dx\right)^{1\over p}\\
&=& {1\over t}\left(\sup_{x\in M} E\left[\int_0^t <XdB_s, T_xF_s(v)>
\right]^q\right)^{1\over q}
|f|_{L^p}.
\end{eqnarray*}
But \begin{eqnarray*}
E\left[\int_0^t<XdB_s, T_xF_s(v)>\right]^q
& \le&  c_p E\left(\int_0^t |T_xF_s(v)|^2 ds\right)^{q/2}
\end{eqnarray*}
by Burkholder-Davis-Gundy's  inequality. Here $c_p$ is a constant.
 So again we have  (\ref{eq: Sobolev 100}).

From  (\ref{eq: manifold 100}) and  corollary \ref{co: 3.2}  we see that
(\ref{eq: formula on manifold}) holds almost everywhere for $f\in L^\infty$ as 
therefore does  (\ref{eq: Sobolev 200}).

\bigskip

\noindent{\bf Example:} Left invariant systems on Lie groups:

Let $G$ be a connected Lie group with identity element ${\bf 1}$ and with
$L_g$ and $R_g$ denoting left and right translation by $G$. Consider a left 
invariant s.d.e.
\begin{equation}
dx_t=X(x_t)\circ dB_t +A(x_t)dt
\label{eq: Lie group}
\end{equation}
with solution $\{g_t: t\ge 0\}$ from ${\bf 1}$. Then (\ref{eq: Lie group})
has solution flow 
$$F_t(u)=R_{g_t}u, \hskip 20pt t\ge 0, u\in G.$$
Take a left invariant Riemannian metric on $G$. Then  by
 (\ref{eq: formula on manifold})  for $f\in BC^1(G)$,  $v_0\in T_{\bf 1}G$, if
 (\ref{eq: Lie group}) is nondegenerate
with $X({\bf 1})\colon R^m\to T_{{\bf 1}}G$ an isometry
\begin{eqnarray*}
dP_tf(v_0)&=&{1\over t} E\left\{f(g_t)\int_0^t <T_{{\bf 1}}R_{{g_s}}(v_0),
X(g_s)dB_s>\right\}\\
&=& {1\over t}E\left\{f(g_t)\int_0^t <ad(g_s)^{-1} (v_0), d\tilde B_s>
_{{\bf 1}}\right\}
\end{eqnarray*}
where $\tilde B_s=X({\bf 1})B_s$. This gives:
$$\nabla \log p_t({\bf 1},y)=
{1\over t} E\left\{ \int_0^t ad(g_s^{-1})^*d\tilde B_s| g_t=y\right\}.$$

\section{For 1-forms}

Let $M$ be a complete Riemannian manifold and $h\colon M\to R$  a smooth
 function with  $L_{\nabla h}$ the Lie 
derivative  in the  direction of $\nabla h$. Let $\triangle^h=\colon \triangle
 +2L_{\nabla h}$ be the Bismut-Witten-Laplacian, and $\triangle^{h,q}$ its
restriction to q-forms. It is then  an essentially
 self-adjoint linear operator on $L^2(M, e^{2h(x)}dx)$  (see \cite{Li.thesis}, 
extending \cite{CH73} from the case $h=0$).  We shall still use
 $\triangle^h$ for its closure and use
$D(\triangle^h)$ for its domain.  By the  spectral theorem, there  is a  smooth
 semigroup $\heatsemi$  solving  the heat equation:

$${\partial P_t\over \partial t}=\half \triangle^h P_t.$$

\bigskip

A stochastic dynamical system (1) is called an {\it h-Brownian system} if it
 has generator $\half \triangle ^h$. Its solution is called an
 {\it h-Brownian  motion}. 

\bigskip

For clarity, we sometime use $P_t^{h,q}$  for 
the  restriction of the semigroup $P_t^h :=\heatsemi$ to $q$-forms.
Denoting  exterior differentiation  by $d$ with suitable domain,  let
$\delta^h$ be the adjoint of $d$ in $L^2(M, e^{2h(x)}dx)$. Then $\triangle^h=
-(d\delta^h +\delta^hd)$, and 
for  $\phi\in D(\triangle^h)$,
\begin{equation} \label{eq: heat}
d(P_t^{h,q}\phi)=P_t^{h, q+1}(d\phi).\end{equation}

\bigskip

Define:
\begin{equation}\label{eq: Bismut 2}
\int_0^t\phi\circ dx_s
=\int_0^t\phi(X(x_s) dB_s) -\half \int_0^t\delta^h\phi(x_s)\, ds,
\end{equation}
for a 1-form $\phi$.
 Theorem \ref{th: Elworthy-Li} has a generalization to closed differential 
forms. It is given in terms of the line integral $\int_0^t\phi\circ dx_s$ 
and a martingale; for it we shall need the following It\^o's formula from
 \cite{ELflow}:

\begin{lemma}[It\^o's  formula for one forms] 
Let $T$ be a stopping time with $T<\xi$, then 
\begin{eqnarray*}
\phi (v_{t\wedge T})
&=&\phi(v_0)+\int_0^{t\wedge T} \nabla \phi\left (X(x_s) dB_s\right)(v_s) 
+ \int_0^{t\wedge T}\phi\left(\nabla X(v_s)  dB_s \right) \\
  &+&\half\int_0^{t\wedge T}\triangle^h\phi(x_s)(v_s)\, ds \\
&+&\half\int_0^t  {\rm trace}\left( \nabla\phi(X(x_s)(-))\nabla X(v_s)(-)
\right)ds. \hskip 10pt \rule{3mm}{3mm}
\end{eqnarray*}

\end{lemma}

\begin{theorem} \label{pr: Bismut 1}
Consider an h-Brownian system. Assume there is no explosion, and
$$\int_0^t E|T_xF_s|^2 ds<\infty, \hskip 6pt \hbox{for each $x$ in $M$.} $$
Let $\phi$ be a closed 1-form    in $D(\triangle^h)\cap L^\infty$,such that
$$\delta P_t\phi=\heatsemif \phi.$$

Then

\begin{equation}
\label{eq: Bismut for 1-forms}
P_t^{h,1}\phi(v_0)={1\over t} E\int_0^t \phi\circ dx_s\int_0^t
<X(x_s)dB_s, TF_s(v_0)>
\end{equation}

\noindent
 for all $v_0\in T_{x_0}M$.
\end{theorem}

\noindent{\bf Proof:}  
Following the proof for a compact manifold as in \cite{ELflow}, let
\begin{equation}
Q_t(\phi)=-\half \int_0^t P_s^{h}(\delta^h \phi)ds. 
\label{eq: Bismut function1}
\end{equation}

\noindent
Differentiate  equation (\ref{eq: Bismut function1}) to get:

$${\partial \over \partial t}Q_t\phi=-\half P_t^{h}(\delta^h\phi).$$
We also have:

\begin{eqnarray*}
d(Q_t\phi)&=&-\half \int_0^t d\delta^h (P_s^h\phi) ds\\
&=&\half \int_0^t \triangle^h(P_s^h \phi)ds\\
&=&P_t^{h}\phi -\phi
\end{eqnarray*}

\noindent
since $d\delta^h(P_s^h\phi)=P_s^{h}(d\delta^h \phi)$ is uniformly continuous 
in $s$ and  $$d(P_s^{h} \phi)=P_s^{h} d\phi=0.$$
  
Consequently:

$$\triangle^h(Q_t(\phi))=-P_t^h(\delta^h\phi)+\delta^h\phi.$$

\noindent
 Apply It\^o's  formula to $(t,x)\mapsto Q_{T-t}\phi(x)$, which is sufficiently
 smooth because $P_s^h\phi$ is, to get:

\begin{eqnarray*}
Q_{T-t}\phi(x_t)&=&Q_T\phi(x_0)+\int_0^t d(Q_{T-s}\phi)(X(x_s)dB_s)\\
&  & +\half \int_0^t \triangle^{h} Q_{T-s}\phi(x_s)ds +\int_0^t{\partial
 \over \partial s}Q_{T-s}\phi(x_s)ds\\
&=& Q_T\phi(x_0)+\int_0^t P_{T-s}^{h} (\phi) (X(x_s)dB_s)-\int_0^t\phi\circ dx_s.
\end{eqnarray*}

\noindent
Setting $t=T$, we obtain:

$$\int_0^T\phi\circ dx_s =Q_T(\phi)(x_0)+\int_0^TP_{T-s}^{h}(\phi)(X(x_s)dB_s),$$

\noindent
and thus 

$$E\int_0^T\phi\circ dx_s\int_0^T<X(x_s)dB_s, TF_s(v_0)>=E\int_0^T P_{T-s}^{h} \phi(TF_s(v_0))ds.$$

\noindent
But

\begin{equation}
E\int_0^T P_{T-s}^{h} \phi(TF_s(v_0))ds
=\int_0^T EP_{T-s}^{h} \phi(TF_s(v_0))ds,
\label{eq: 1-form 2}
\end{equation}

\noindent
 by Fubini's theorem, since

\begin{eqnarray*}
\int_0^T E|P_{T-s}^{h} \phi (TF_s(v_0))| ds
&\le& |\phi|_\infty \int_0^T E|TF_T(v_0)| ds<\infty.
\end{eqnarray*}

\noindent
Next notice:

$$E P_{T-s}^{h}\phi (  TF_s(v_0)  ) =E\phi(TF_T(v_0))=P_T^{h}\phi(v_0)
$$

\noindent
from the strong Markov property. We get: 
$$P_T^{h,1}\phi(v_0)={1\over T} E\left\{\int_0^T\phi\circ dx_s\int_0^T
<XdB_s, TF_s(v_0)>\right\}.$$
\hfill \rule{3mm}{3mm}

\bigskip
\noindent{\bf Remark:}\label{remark}
If we assume  $\sup_x E|T_xF_t|^2<\infty$ for each $t$ the result holds for
all $\phi\in D(\triangle^h)$: first   we have
 $\delta P_t\phi=\heatsemi \phi$ for $\phi\in L^2$ by continuity and  also
 equation (\ref{eq: 1-form 2})  holds from the following argument:
 
\begin{eqnarray*}
&&\int_0^T E |P_{T-s}^{h}\phi (TF_s(v))| ds
\le\int_0^T E|\phi(TF_T(v))| \, ds\\
&&\le  E|T_xF_T(v)|^2 \sup_x  \left(\int_0^T E|\phi|_{F_T(x)}^2 ds\right) \\
\end{eqnarray*}

\noindent
But $\int_ME|\phi|_{F_T(x)}^2 e^{2h} dx =\int |\phi|^2 e^{2h} dx<\infty$. So 
$E|\phi|^2_{F_T(x)}<\infty$ for each $x$ by the continuity of
 $E|\phi|_{F_T(x)}^2=P_T(|\phi|^2)(x)$ in $x$.

\bigskip

\begin{corollary}
Suppose $|\nabla X|$ is bounded and for all $v\in T_xM$, all $x\in M$

\noindent
$ Hess(h)(v,v)-\half Ric_x(v,v)\le c|v|^2$ for some constant $c$. Then
(\ref{eq: Bismut for 1-forms}) holds
for all closed 1-forms in $D(\triangle^h)$.  
\end{corollary}

\noindent{\bf  Proof:} By  lemma A\ref{le: 7.2}, we have
$\sup_x E|T_xF_t|^2<\infty$ and $E\sup_{s\le t}|T_xF_s|<\infty$.  Thus
proposition A\ref{pr: 7.6}  shows that
 $P_t^{h,1} \phi=\delta P_t(\phi)$.
Theorem \ref{pr: Bismut 1} now applies. \hfill\rule{3mm}{3mm}

\bigskip

\noindent{Remark:}
Note that if $\phi=df$, formula (\ref{eq: Bismut for 1-forms}) reduces to
(\ref{eq: formula on manifold}) using (\ref{eq: heat}).

\section{The Hessian flow}

{\bf A.}  Let $Z=A^X$   as in section 3.   Let $x_0\in M$ with 
 $\{x_t: 0\le t<\xi\}$  the solution to (1) with initial value $x_0$  and
explosion time $\xi$.
Let $W_t^Z$ be the solution flow to  the covariant differential
 equation along $\{x_t\}$:
\begin{equation}
{ DW_t^Z(v_0) \over \partial t}
=-\half {\rm Ric}^\#(W_t^Z(v_0),-)+\nabla Z(W_t^Z(v_0))
\end{equation}
with $W_0^Z(v_0)=v_0$.  It is called the Hessian flow. Here Ric denotes the 
Ricci curvature of the manifold, and  $\#$ denotes the relevant raising or
 lowering of indices so that $\hbox{Ric}^\#(v,-)\in T_xM$ if $v\in T_xM$. 
For $x\in M$ set

$$\rho(x)=\inf_{|v|\le 1}\{Ric_x(v,v)-2\nabla Z(x)^\#(v,v)\}.$$

 The following is a generalization of a result in \cite{ELflow}.

\begin{proposition}\cite{Li.thesis}
 Let $Z=\nabla h$ for $h$ a smooth function on $M$.
Suppose for some $T_0>0$,
$$E\sup_{t\le T_0}\chi_{t<\xi(x)}\exp^{-\half \int_0^t \rho(F_s(x))ds}<\infty,
\hskip 6pt 0\le t\le T_0$$
 Then for a closed bounded $C^2$  1-form $\phi$, we have  for $0<t\le T_0$:

\begin{equation}
\label{eq: Bismut for Hessian}
P_t^{h,1}\phi(v_0)={1\over t} E\int_0^t \phi\circ dx_s\int_0^t
<X(x_s)dB_s, W_s^Z(v_0)>.
\end{equation}
\end{proposition}

The proof is as for (\ref{eq: Bismut for 1-forms}) with $TF_t$, just noticing
 that under the   conditions of the proposition, the s.d.e. does not explode
 and   $P_t^{h,1}\phi=E\phi(W_t^h)$ for bounded 1-forms $\phi$ 
(see e.g. \cite{ELflour} and \cite{application}).

\bigskip

\noindent{\bf Remark:}
 Taking $\phi=df$, we obtain, by (\ref{eq: heat}),
\begin{equation}\label{eq: Hessian}
dP_tf(v_0)={1\over t}Ef(x_t)\int_0^t<W_s^Z(v_0), X(x_s)dB_s>\end{equation}
which leads to  Bismut's formula \cite{ELflow} for $\nabla \log p_t(-,y)$ 
(proved there for $Z=0$ and $M$ compact). In fact (\ref{eq: Hessian})
can be proved directly, without assuming $Z$ is a gradient, by our basic
method:
Let $\phi_t=d(P_tf)$, then it solves
${\partial \phi_t\over \partial t}=\half \triangle^1\phi_t +L_{\nabla Z}\phi_t$
 since  $P_tf$ solves  ${\partial g \over \partial t}=\half \triangle g
 +L_{\nabla Z } g$. Then It\^o's formula  (as in \cite{ELflow}) applied  
to $\phi_{t-s}(W_s^Z(v_0))$ shows that $\phi_t(v_0)=E\phi_0(W_t^Z(v_0))$
and our usual method can be used.
Furthermore if $\rho$ is bounded from below so that $|W_t^Z|$ is bounded as
 in \cite{ELflow}, then  (\ref{eq: Hessian})  holds for  bounded measurable
 functions.

\bigskip

Note that it was shown,   in  \cite{EL-YOR},  that for a  gradient  system on
  compact  $M$, $E\{v_t| x_s: 0\le s \le t\}=W_t^h(v_0)$. 
Recall that a {\it gradient system} is given by 
 $X(\cdot)(e)=\nabla< f(\cdot),e>$,
$e\in R^m$ for $f\colon M\to R^m$ an isometric embedding.  
This relation between the derivative flow and the Hessian flow holds 
for noncompact manifolds  if  $E\int_0^t|\nabla X(x_s)|^2 |v_s|^2ds<\infty$.

\bigskip

\noindent {\bf B.}
Let $V_\cdot(\cdot)\colon [0,\infty)\times R^n \to R$ be continuous, $C^1$
in $x$ for each $t$ and bounded above with derivative $dV$ bounded  on each
 $[0,T]\times R^n$. Consider the following equation with potential $V$
$${\partial u_t\over \partial t}=\half \triangle u_t+L_Zu_t+V_tu_t.$$
Assume that the  s.d.e. (1) does not explode. By the corresponding argument to
that used for the case $V\equiv 0$,  we get for $v_0\in T_{x_0}M$

\begin{eqnarray*}
du_t(v_0)&=&Eu_0(x_t)\exp^{\int_0^tV_{t-s}(x_s)ds}
\int_0^t dV_{t-s}(W_s^Z(v_0))ds\\
&+& Edu_0(W_t^Z(v_0))\exp^{\int_0^tV_{t-s}(x_s)ds}
\end{eqnarray*}
provided that $-\half Ric^\#+\nabla Z$ is bounded above as a linear 
operator and $u_0$ is $BC^1$. From this the  analogous proof to that of
 (\ref{eq: Remark 10})  gives:

\begin{theorem}
 Assume nonexplosion and suppose $-\half Ric^\#+\nabla Z$ is
 bounded above and dV is bounded. Then  for $u_0$ bounded measurable and $t>0$,

\begin{equation}\label{eq: Remark 20} \begin{array}{ll}
du_t(v_0)&={1\over t}Eu_0(x_t)\exp^{\int_0^t V_{t-s}(x_s)ds}
\int_0^t<W_s^Z, X(x_s)dB_s>\\
&+{1\over t}E u_0(x_t)\exp^{\int_0^t V_{t-s}(x_s)ds}\int_0^t
(t-r) dV_{t-r}(x_r)(W_r^Z)dr.
\end{array}\end{equation}

\end{theorem}

\section{For higher order forms and gradient Brownian systems}

Recall that a gradient h-Brownian system is a gradient system with  
 $A(x)=:\nabla  h(x)$.   For such systems
$\sum_1^m \nabla X^i(X^i)=0$.  We shall assume there is no explosion 
 as before.

If $A$ is a linear map from  a vector space $E$ to $E$,  then 
$(d\Lambda)^q A$ is the map from $E\times \dots \times E$ to $E\times \dots 
\times E$ defined as  follows:
$$(d\Lambda)^q A(v^1,\dots, v^q)
=\sum_{j=1}^q (v^1, \dots, Av^j, \dots, v^q).$$

\noindent
 Let $v_0=(v^1_0, \dots, v_0^q)$, for
 $v_0^i\in T_{x_0}M$.  Denote by $v_t$ the $q$ vector induced by $TF_t$: 
 $$v_t=(TF_t(v_0^1), TF_t(v_0^2),\dots, TF_t(v_0^q)).$$

\begin{lemma}\cite{ELflow}
Let  $\theta$ be a $q$ form. Then,  for a gradient h-Brownian system,

  \begin{eqnarray*}
\label{eq: Ito formula for gradient}
\theta(v_t)&=&\theta(v_0)+\int_0^t\nabla \theta(X(x_s)dB_s)(v_s)\\
&+&\int_0^t\theta\left((d\Lambda) ^q (\nabla X(-)dB_s)(v_s)\right)
+\int_0^t \frac 1 2 \triangle ^{h,q}(\theta)(v_s)ds. \hfill \rule{3mm}{3mm}
\end{eqnarray*}
\end{lemma}

\noindent  Recall that  if $\theta$ is a $q$ form, then

\begin{equation}
 (\delta P_t)\theta(v_0)\index{\delta P_t}=E\theta(v_t)
\end{equation}
 where defined.   Define a (q-1) form $\int_0^t \theta\circ dx_s$ by

\begin{equation}\label{eq: definition of path integral}\begin{array}{ll}
\int_0^t\theta \circ d x_s(\alpha_0)=:& 	
	{1\over q}\int_0^t\theta\left(X(x_s)dB_s, TF_s(\alpha_0^1), \dots,
TF_s(\alpha_0^{q-1})\right)\\
	 &-{1\over 2}\int_0^t\delta^h \theta \left(TF_s(\alpha_0^1), \dots, 
TF_s(\alpha_0^{q-1})\right)\, ds
\label{eq: h-forms 2}
\end{array}\end{equation}
for $\alpha_0=(\alpha_0^1, \dots, \alpha_0^{q-1})$ a (q-1) vector.
Then we have the following extension of theorem \ref{pr: Bismut 1}:

\begin{theorem} Let $M$ be a complete Riemannian manifold. Consider a gradient
h-Brownian system on it. Suppose it has no explosion and for each $t>0$ and
$x\in M$,
$$\int_0^t  E|T_xF_s|^{2q} ds <\infty.$$
Let $\theta$ be a closed  bounded  $C^2$   $q$ form in $ D(\triangle^{h,q})$ 
 with $$\delta P_t\theta= P_t^{h,q} \theta.$$
 Then:

\begin{equation}\label{eq: Bismut higher forms}
 (P_t^{h,q}\theta)_{x_0}=  \frac 1t
 E \int_0^t<X(x_s)dB_s, T_{x_0}F_s(\cdot)> \wedge \int_0^t\theta\circ dx_s.
\end{equation}
\end{theorem}

\noindent{\bf Proof:}  
  Let $Q_t\theta$ be the (q-1) form given by
\begin{equation}
Q_t(\theta)(\alpha_0)= -\frac12\int_0^t(\delta^{h}P_s^{h,q}\theta)(\alpha_0)ds,
\end{equation}
for $\alpha_0\in \wedge^{q-1}T_xM$.

Notice that  $P_t^{h,q}(\theta)$ is smooth on $[0,T]\times M$ by parabolic 
regularity, so

\begin{eqnarray*}
\frac \partial {\partial t} Q_t(\theta)&=&-\frac 12 \delta^h(P_t^h\theta),\\
	d(Q_t(\theta))&=&-\frac 12\int_0^td\delta^h(P_s^{h,q}\theta)ds, \\
	\delta^h Q_t(\theta)&=&-\half \int_0^t \delta^h\delta^h (P_s^{h,q}\theta)\, ds
=0.
\end{eqnarray*}

\noindent
 In particular, 
 \begin{equation}\label{eq: cohomology}
d(Q_t(\theta))=  \frac 12\int_0^t\triangle^{h,q}(P_s^{h,q}\theta)ds       
	=P_t^{h,q}\theta - \theta
\end{equation}

\noindent
since $\triangle^{h,q} \theta= -d\delta^h \theta$.
Therefore:

  $$\triangle ^{h,q-1} (Q_t(\theta))
=-P_t^{h, q-1}(\delta^h\theta) +\delta^h\theta.$$
 
Next we apply It\^o's formula (the previous lemma)  to $(t,\alpha)\mapsto
  Q_{T-t}(\theta)(\alpha)$, writing
 $\alpha_t=\left(TF_t(\alpha_0^1),\dots, TF_t(\alpha_0^{q-1})\right)$:

\[\begin{array}{ll}
	Q_{T-t}\theta(\alpha_t)=& Q_T\theta(\alpha_0)
	+\int_0^t\nabla Q_{T-s}\theta(X(x_s)dB_s)(\alpha_s)    \\
	&+      \int_0^tQ_{T-s}\theta\left(
(d\wedge)^{q-1}(\nabla X(-)dB_s)(\alpha_s)  \right)   \\
	&+\frac 12\int_0^t\triangle^h Q_{T-s}\theta(\alpha_s)ds
	+\int_0^t\frac \partial{\partial s}(Q_{T-s})\theta(\alpha_s)ds.
\end{array}\]

From  the  calculations above  we get:  

\begin{eqnarray*}
	Q_{T-t}\theta(\alpha_t)&=&Q_T\theta(\alpha_0)
	+\int_0^t\nabla Q_{T-s}\theta(X(x_s)dB_s)(\alpha_s)     \\
	&+&  \int_0^tQ_{T-s}\theta\left(
(d\wedge)^{q-1}(\nabla X(-)dB_s)(\alpha_s)\right)   \\
	&+&   \frac 12 \int_0^t \delta ^h \theta(\alpha_s) ds.
\end{eqnarray*}

 By definition and the equality above,

\begin{equation}\begin{array}{ll}
\int_0^T\theta\circ dx_s(\alpha_0)=&Q_T\theta(\alpha_0)+ 
\frac 1q\int_0^T\theta(X(x_s)dB_s, \alpha_s) \\
	 &+ \int_0^T\nabla Q_{T-s}\theta(X(x_s)dB_s)(\alpha_s)    \\ 
	 &+\int_0^T Q_{T-s}\theta\left(
(d\wedge)^{q-1}(\nabla X(-)dB_s)(\alpha_s)  \right).     
\label{eq: formulae 2}
\end{array}\end{equation}

We will calculate the expectation of each term of $\int_0^t \theta\circ dx_s$
in (\ref{eq: formulae 2})   after wedging with
 $\int_0^T<X(x_s)dB_s, TF_s(-)>ds$.  The first term clearly  vanishes.
 The last term vanishes as well for a gradient h-Brownian
 system since $\sum_i\nabla X^i(X^i(-))=0$. 

Take $v_0=(v_0^1, \dots, v_0^q)$.  Write $ v_s^i=TF_s(v_0^i)$, and
denote by  $w_s(\cdot)$ the linear map:

 $$w_s(\cdot)=\overbrace{ (TF_s(\cdot), \dots, TF_s(\cdot))}^{q-1}.$$ 
Then
\begin{eqnarray*}
	&&{1\over q}  E\int_0^T\theta(X(x_s)dB_s, w_s(\cdot)) 
		\wedge\int_0^T<X(x_s)dB_s, TF_s(\cdot)>(v_0)        \\
&=&{1\over q} \sum_{i=1}^q (-1)^{q-i} E\int_0^T\theta(v_s^i, v_s^1, \dots, 
\widehat{v_s^i} , \dots , v_s^q) \, ds  \\
&=&{1\over q} \sum_{i=1}^q (-1)^{q-i}(-1)^{i-1} E\int_0^T \theta(v_s^1,\dots , v_s^q)\,ds  \\
	&=& (-1)^{q-1} E\int_0^T\theta(v_s^1, \dots, v_s^q)ds       \\
	&=&(-1)^{q-1} \int_0^T P_s^{h} \theta (v)ds. 
\end{eqnarray*}

\noindent
The last step uses the assumption: $\int_0^T  E|T_xF_s|^{2q} ds <\infty$. 
Similar calculations show:

\begin{eqnarray*}
	 &&E\{\int_0^T\nabla Q_{T-s}\theta(X(x_s)dB_s)(w_s(\cdot))
		\wedge \int_0^T<X(x_s)dB_s, TF_s(\cdot)>\}(v_0)       \\
&=& \sum_{i=1}^q (-1)^{q-i} E\int_0^T \nabla (Q_{T-s}\theta)(v_s^i)
		 (v_s^1 , \dots  \hat{v_s^i} , \dots, v_s^q) \, ds     \\
	&=& (-1)^{q-1}E\int_0^T (d(Q_{T-s}\theta))(v_s^1, \dots, v_s^q)\, ds, \\
&=&(-1)^{q-1} \int_0^T P_s^{h} \left(P_{T-s}^h(\theta)-\theta\right)\, (v)\, ds  \\
	&=& (-1)^{q-1} \left[T(P_T^{h}\theta)(v)  - \int_0^T P_s^h\theta(v)ds\right].
\end{eqnarray*}

\noindent
Comparing these with $(\ref{eq: formulae 2})$, we have: 

$$ P_T^{h,q}\theta ={1\over T} E \int_0^T<X(x_s)dB_s, TF_s(\cdot)>
\wedge\int_0^T \theta\circ dx_s. $$
 \hfill \rule{3mm}{3mm}
\noindent

Note: With an additional condition: $\sup_{x\in M}E|T_xF_s|^{2q}<\infty$,
the formula in the above proposition holds for forms which are not necessarily 
bounded. See the remark at the end of section 4.
\bigskip

Recall that $\rho(x)$ is the distance function between $x$ and a fixed point in
 $M$, and  ${\partial h\over \partial \rho}:=dh(\nabla \rho)$.

\begin{corollary}
 Consider a gradient h-Brownian system. Formula (\ref{eq: Bismut higher forms})
  holds  for a  closed  $C^2$ q-form in $D(\triangle^h)$, 
if one of the following conditions holds:
\begin{enumerate}\item
The  related second fundamental form is bounded and $\half Ric -Hess(h)$ is 
bounded from below;
\item
The second fundamental form is bounded by $c[1+\ln(1+\rho(x))]^{\half}$,
and also  $${\partial h\over \partial \rho}\le c[1+\rho(x)],$$
$$ Hess(h)(x)(v,v)\le c[1+\ln(1+\rho(x))]|v|^2.$$
\end{enumerate}
\end{corollary}

\noindent {\bf Proof:}
This follows since,  \cite{ELflow}, for $v_1, v_2\in T_xM$, and $e\in R^m$ 
$$<\nabla X(v_1)e, v_2>_x=<\alpha(v_1,v_2),e>_{R^m}.$$
 Lemma A\ref{le: 7.2} and lemma A\ref{new 7.3}
  give $E\sup_{s\le t}|T_xF_s|^{2q}<\infty$ for all $q$.
 The second part of proposition A\ref{pr: 7.6}   now gives
$\delta P_t\theta= P_t^{h,q} \theta$. So the   conditions of the  theorem
 are satisfied, with the remark above used to avoid assuming $\theta$ is
bounded.  \hfill\rule{3mm}{3mm}

\bigskip

We now have the extension  of our basic differentiation result to the case of
   q-forms.

\begin{corollary}
 Consider a gradient h-Brownian system  on a complete
 Riemannian manifold. Suppose there is no explosion  and 
$\int_0^t E|T_xF_s|^{2q}ds<\infty$. Let $\phi$  be  a q-1 form  such that
$d\phi$ is a bounded $C^2$ form in $D(\triangle^{h,q})$ with
$P_t^{h,q}(d \phi)=\delta P_t(d\phi)$. Then

\begin{equation}\label{differentiation for forms}
 d\left(P_t^{h,q-1}(\phi)\right)=\frac 1t
E\left( \int_0^t<X(x_s)dB_s, TF_s(\cdot)>\wedge 
\phi(\overbrace{TF_t(\cdot),\dots TF_t(\cdot)}^{q-1})\right).\\
\end{equation}
\end{corollary}

\noindent   {\bf Proof:}  
  By (\ref{eq: definition of path integral}), if $\theta=d\phi$
 
\begin{equation}\label{eq: Higher 19}\begin{array}{ll}
\int_0^t\theta \circ dx_s(-) =&\frac 1q \int_0^t d\phi(X(x_s)dB_s,
 \overbrace{TF_s(\cdot), \dots, TF_s(\cdot)}^{q-1})(-)\\
&+\frac 12 \int_0^t \triangle^h\phi(\overbrace{TF_s(\cdot), \dots, TF_s(\cdot}^{q-1})	 )(-)ds.
 \end{array}\end{equation}

On the other hand, if $\alpha_0=(\alpha_0^1, \dots, \alpha_0^{q-1})$ for
 $\alpha_0^i\in T_{x_0}M$, then  by It\^o's formula,

\begin{equation}
\begin{array}{ll}
\phi\left(TF_t(\alpha_0^1), \dots, TF_t(\alpha_0^{q-1})\right)
=&\phi(\alpha_0)
+\int_0^t \nabla \phi(X(x_s)dB_s)
\left(TF_s(\alpha_0^1), \dots, TF_s(\alpha_0^{q-1})\right) \\
&+\half\int_0^t\triangle^h\phi\left(TF_s(\alpha_0^1), \dots,
 TF_s(\alpha_0^{q-1}) \right)ds.\\
\end{array}\label{eq: Higher 30}\end{equation}

However

\begin{eqnarray*}
	&E&   \int_0^t <X(x_s)dB_s,\, TF_s(\cdot)>\wedge
\int_0^t d\phi\left(X(x_s)dB_s, \, TF_s(\cdot),\dots, TF_s(\cdot)\right )  \\
	&=& q E \int_0^t <X(x_s)dB_s, TF_s(\cdot)>\wedge
\int_0^t\nabla\phi (X(x_s)dB_s)\left( TF_s(\cdot),\dots, 
TF_s(\cdot)\right).
\end{eqnarray*}

Compare (\ref{eq: Higher 19}) and equation (\ref{eq: Higher 30}) to obtain:
\begin{eqnarray*}
&&E\int_0^t  <XdB_s, TF_s(-)>\wedge \int_0^t d\phi\circ dx_s\\
&=&E\int_0^t   <XdB_s, TF_s(-)>\wedge
\int_0^t \nabla\phi\left (XdB_s\right)\overbrace{
\left(TF_s(-),\dots,  TF_s(-)\right)}^{q-1}\\
&+&\half\int_0^t   <XdB_s, TF_s(-)>\wedge
\int_0^t\triangle^h \phi(\overbrace{TF_t(-),\dots TF_t(-)}^{q-1})ds\\
&=&E\left(\int_0^t<XdB_s, TF_s(-)>\wedge
 \phi(\overbrace{TF_t(-),\dots TF_t(-)}^{q-1})\right). \end{eqnarray*}
This gives the required result by the formula for $P_t^{h,q}(d\phi)$ in the
 previous theorem.
 \hfill \rule{3mm}{3mm}

\bigskip

\noindent{\bf Remarks:}

(i) This can be proved directly as for the case $q=0$
 in theorem \ref{th: Elworthy-Li}.

(ii) Equation (\ref{differentiation for forms})
can be given the following interpretation:

\noindent
Our stochastic differential equation determines a 1-form valued process
$\Psi_{t}=\Psi_{t}^{X,A}, t\ge 0$ given by
$$\Psi_{t,x_0}(v_0)=\int_0^t <X(x_s)dB_s, T_{x_0}F_s(v_0)>$$
i.e. 
$$\Psi_{t, x_0} =\int_0^t (T_{x_0}F_s)^*\left(<X(x_s)dB_s, ->_{x_s}\right)$$
(so for each $x_0$, $\{\Psi_{t,x_0}: t\ge 0\}$ determines a local martingale
on $T_{x_0}^*M$ with tensor quadratic variation given by 
$\int_0^t T_\cdot F_s^*T_\cdot F_sds$. Note that the Malliavin covariance
 matrix is given by 
 $\int_0^t \left(T_\cdot F_s^*T_\cdot F_s\right)^{-1}ds \hskip 4pt$ ). In fact
$\Psi_t$ is exact: $\Psi_t=d\psi_t$ where $\psi_t: M\times \Omega\to R$ is
 given by
$$\psi_t(x)=\int_0^t <f\left(F_s(x)\right), dB_s>_{R^m}$$
for $f: M\to R^m$ the given embedding. 

 Equation (\ref{differentiation for forms}) states

\begin{eqnarray*}
dP_t^{h,q-1}\phi&=&{1\over t}E\{ \Psi_t\wedge  (F_t)^*\phi\}\\
&=&{1\over t}E\{ d\psi_t\wedge  (F_t)^*\phi\}\\
\end{eqnarray*}

(iii) Note that (\ref{eq: cohomology}) gives an explicit cohomology between
$P_t^{h,q}\theta$ and $\theta$.

\section*{Appendix: Differentiation under the expectation}

Consider the  stochastic differential equation:
\begin{equation}\label{eq: 0}
dx_t=X(x_t)\circ dB_t+A(x_t)dt
\end{equation}
on a complete n-dimensional Riemannian manifold. We need the following result
 on the existence of a partial flow  taken from \cite{ELbook}, following
Kunita:

\begin{theorem*} $\label{th: partial flow}$
Suppose $X$, and $A$ are  $C^{r}$, for $r\ge 2$. Then there is a partially
 defined flow $(F_t(\cdot),\xi(\cdot))$ such that for each $x\in M$, 
$(F_t(x), \xi(x))$  is a maximal solution  to $(\ref{eq: 0})$ with lifetime
$\xi(x)$ and if 

$$M_t(\omega)=\{x\in M, t<\xi(x,\omega)\},$$
then there is a set $\Omega_0$ of full measure such that for all $\omega\in 
\Omega_0$:
\begin{enumerate}
\item
$M_t(\omega)$ is open in $M$ for each $t>0$, i.e. $\xi(\cdot,\omega)$ is lower
 semicontinuous.
\item
$F_t(\cdot,\omega): M_t(\omega)\to M$ is in $C^{r-1}$ and is a diffeomorphism
 onto an open subset of $M$. Moreover the map : $t\mapsto F_t(\cdot,\omega)$ is
continuous into $C^{r-1}(M_t(\omega))$, with the topology of uniform 
convergence on compacta of  the first r-1 derivatives.
\item
Let $K$ be a compact set and $\xi^K=\inf_{x\in K} \xi(x)$. Then

\begin{equation}
\lim_{t\nearrow \xi^K(\omega)} \sup_{x\in K} d(x_0, F_t(x))=\infty
\end{equation}

\noindent
almost surely on the set $\{\xi^K<\infty\}$. (Here $x_0$ is a fixed point of
 $M$ and $d$ is any complete  metric on $M$.)

\end{enumerate}
\end{theorem*}

From now on, we shall use $(F_t,\xi)$ for the partial flow defined in 
theorem A\ref{th: partial flow} unless otherwise stated. 

Recall that a stochastic differential equation is called {\it strongly 
p-complete} if its solution can be chosen to be jointly continuous in time and 
space for all time when  restricted to a smooth singular p-simplex. A singular 
p-simplex is a map $\sigma$ from a standard p-simplex to $M$. We also use the 
term 'singular p-simplex' for its image. If a s.d.e. is strongly p-complete, 
$\xi^K=\infty$ almost surely for each smooth singular p-simplex $K$
\cite{flow}. 

Let $x\in M$, and $v\in T_xM$. Define $H_p$ as follows:
\begin{eqnarray*}
H_p(x)(v,v)&=&2<\nabla A(x)(v),v>+\sum_1^m <\nabla^2 X^i(X^i,v),v>
+ \sum_1^m|\nabla X^i(v)|^2\\
&+&  \sum_1^m<\nabla X^i\left(  \nabla X^i(v)  \right),v>+
(p-2)\sum_1^m{1\over |v|^2} <\nabla X^i(v),v>^2.
\end{eqnarray*}

 There are simplifications of $H_p$:

\noindent
For  s.d.e. (3) on $R^n$,
\begin{eqnarray*}
H_p(x)(v,v)= 2<DZ(x)(v),v>+ \sum_1^m|D X^i(v)|^2+
  (p-2)\sum_1^m{ <D X^i(v),v>^2\over |v|^2}.
\end{eqnarray*}

\noindent
For (1) with generator $\half \triangle +L_{Z}$,

\begin{eqnarray*}
H_p(x)(v,v)&=&-Ric_x(v,v)+ 2<\nabla Z(x)(v),v>+
 \sum_1^m|\nabla X^i(v)|^2\\
&+&(p-2)\sum_1^m{1\over |v|^2} <\nabla X^i(v),v>^2.
\end{eqnarray*}

There are the following lemmas from \cite{flow}:

\begin{lemma*}\label{le: 7.2}
Assume the stochastic differential equation (1) is complete. Then

(i).  It is strongly 1-complete if  $H_1(v,v)\le c|v|^2$ for
 some constant $c$. Furthermore if also  $|\nabla X|$ is bounded, then it is
strongly complete and 

\noindent
 $\sup_x E\left(\sup_{s\le t}|T_xF_s|^p\right)$ is finite for all $p>0$ 
  and  $t>0$.

(ii) Suppose $H_p(v,v)\le c|v|^2$, then
$\sup_{x\in M}E|T_xF_t|^p \le k\exp^{cp{t\over 2}}$ for $t>0$. Here $k$ is a
 constant independent of $p$.
\end{lemma*}

For a more refined result, let $c$ and $c_1$ be two constants, let
 $ \rho(x)$ be the distance between $x$
 and a fixed point $p$ of $M$, and assume $\A=\half \triangle +Z$.

\begin{lemma*}\cite{flow}\label{new 7.3}
Assume that the Ricci curvature at each point $x$ of $M$ is bounded from below
 by $-c(1+\rho^2(x))$.
Suppose $dr(Z(x))\le c[1+\rho(x)]$, then there is no explosion.
If furthermore  $|\nabla X(x)|^2\le c\left[1+\ln(1+\rho(x))\right]$, and
$$Ric_x(v,v)-2<\nabla Z(x)(v),v>\ge c_1[1+\ln(1+\rho(x))]|v|^2,$$
then the system is strongly complete and
$$\sup_{x\in K}E\left(\sup_{s\le t}|T_xF_s|^p\right)<k_1\exp^{k_2t}.$$
for all compact sets $K$. Here $k_1$ and $k_2$ are constants independent of
 $t$.
\end{lemma*}

We will first use strong 1-completeness to differentiate under expectations
in the sense of distribution. For this furnish $M$ with a complete Riemannian
 metric and let $dx$ denote the 
corresponding volume measure of $M$. Let $\Lambda$ be a smooth vector field on 
$M$. For $f\in L_{loc}^1(M,R)$, the space of locally integrable functions on 
$M$,  we say that $g\in L_{\rm loc}^1(M,R)$ is the {\it weak Lie derivative}
 of $f$  in the direction $\Lambda$ and write
$$g=\L_\Lambda f, \hskip 20pt \hbox{weakly}$$
if for all $\phi:M\to R$ in $C_K^\infty$, the space of smooth functions 
with compact support,  we have:
$$\int_M\phi(x)g(x)dx
=-\int_Mf(x)\left[<\nabla \phi(x), \Lambda(x)>_x+\phi(x){\rm div}\Lambda(x)
\right]dx.$$
A locally integrable 1-form $\psi$ on $M$ is the {\it weak derivative} of $f$
$$df=\psi \hskip 20pt {\rm weakly}$$
if $\psi(\Lambda(\cdot))=\L_\Lambda f$ weakly for all $C_K^\infty$ vector
 fields  $\Lambda$ on $M$.

Let $\Lambda$ be a $C_K^\infty$ vector field on $M$ and for each $x$ in $M$
let $K(x)$ be the integral curve of $\Lambda$ through $x$.

\begin{lemma*}\label{le: 7.4}
Suppose the s.d.e. (\ref{eq: 0}) is complete. Then for $t\ge 0$,

(i) With probability one $M_t(\omega)=\{x: t<\xi(x,\omega)\}$ has full
measure in $M$. In particular $f\circ F_t(-,\omega)$ determines an element of
$L_{\rm loc}^1(M,R)$ with probability one for each bounded measurable
$f: M\to R$.

\bigskip

If also (\ref{eq: 0}) is strongly 1-complete and $f$ is $BC^1$ then with probability 1:

(ii)  $t<\xi^K(\omega)$ for any compact  subset $K$ of $K(x)$  for
 almost all $x$ in $M$;

\noindent
and

(iii) the Lie derivative $\L_\Lambda(f\circ F_t(-,\omega))$ exists almost
everywhere on $M$ in the classical sense, is equal to the Lie derivative
in the weak sense almost everywhere, and
$$\L_\Lambda(f\circ F_t(-,\omega))=df\circ TF_t(-,\omega)
(\Lambda(-)) \hskip 16pt {\rm weakly}.$$
\end{lemma*}

\noindent
{\bf Proof:} Completeness of (\ref{eq: 0}) implies that $\xi(x,\omega)=\infty$
with probability 1 for each $x$ in $M$ so that
$\{(x,\omega)\in M\times \Omega: t<\xi(x,\omega)\}$ has full $\lambda
\otimes P$ measure. Fubini's theorem gives (i). The same argument applied to
$\{(x,\omega)\in M\times \Omega: t<\xi^{K(x)}(\omega)\}$ yields (ii).

From (ii) we know that if $f\in BC^1$, then $f\circ F_t(-,\omega)$ is $C^1$ on
 almost all $\{K(x): x\in M\}$ with probability one. In particular it is
 absolutely  continuous along the trajectories of $\Lambda$ through almost all
 points of $M$ with probability one. It follows e.g. by Schwartz
 \cite{Schwartz66}  chapter 2 section 5 that 

\noindent
$\L_\Lambda\left( f\circ F_t(-,\omega)\right)$
 exists almost  everywhere. However at each point $x$ of $M_t(\omega)$
 this classical  derivative is just $df\circ T_xF_t(-,\omega)(\Lambda(x))$,
which is in $L^1_{loc}$. By \cite{Schwartz66} it is therefore equal to the 
weak Lie derivative almost everywhere, with  probability 1. 
\hfill\rule{3mm}{3mm}

\begin{theorem*}\label{th: 7.5}
Suppose the stochastic differential equation (\ref{eq: 0}) is strongly 
1-complete and $E|T_xF_t|\in L^1_{loc}$ in $x$.  Then for $f$ in $BC^1$, 
$P_tf$ has weak derivative given by
\begin{equation}\label{eq: appendix 500}
d(P_tf)=\delta P_t(df) \hskip 24pt \hbox{weakly}
\end{equation}
In particular this holds if $H_1(v,v)\le c|v|^2$.
\end{theorem*}

\noindent {\bf Proof:}
Let $\Lambda$ be a $C_K^\infty$ vector field on $M$. Then by lemma 
 A\ref{le: 7.4} and Fubini's theorem:
\begin{eqnarray*}
\int_MP_tf(x){\rm div}\Lambda(x)dx
&=& \int_M Ef(F_t(x)){\rm div}\Lambda(x)dx\\
&=& E\int_M f\circ F_t(x, \omega){\rm div}\Lambda(x)dx\\
&=& -E\int_M \L_\Lambda(f\circ F_t(-, \omega))(x)dx\\
&=& -E\int_M df\circ T_xF_t(-,\omega)(\Lambda(x))dx\\
&=&-\int_M \delta P_t(df) (\Lambda(x)) dx
\end{eqnarray*}
as required. The last part comes from lemma A\ref{le: 7.2}.
\hfill\rule{3mm}{3mm}

\bigskip

\noindent{\bf Remark:}
Under the conditions of the theorem it follows as in \cite{Schwartz66} that
 the derivatives $\L_\Lambda(P_tf)$ exist in the
classical sense a.e. for each smooth vector field $\Lambda$ and are
given by $\delta P_t(df)(\Lambda(\cdot))$.

 If also the stochastic differential equation (\ref{eq: 0}) is non-degenerate
 (so that its  generator is elliptic) and $x\to E|T_xF_t|$ is continuous on
each compact set, then by parabolic regularity and a direct proof in
 \cite{flow} equation  (\ref{eq: appendix 500}) holds
in the classical sense at all points of $x$.

\bigskip

In the elliptic case there are the following criteria:

\begin{proposition*}\label{pr: 7.6}   \cite{ELflour}\cite{Li.thesis}
For a complete h-Brownian system on a  complete Riemannian manifold:

(i)  suppose $E\sup_{s\le t}|T_xF_s|<\infty$ for all $x\in M$ and $t>0$,
 then  for every bounded $C^2$, closed 1-form $\phi_0$,  $\delta P_t(\phi_0)$
is the unique solution  to the heat equation
${\partial \phi_t\over \partial t}=\half \triangle^{h,1}\phi_t$
with initial condition $\phi_0$.
 If $\phi=df$, this gives
$dP_tf(x)=\delta P_t(df)(x)$ for all $x$ and
for all bounded $C^3$ functions with bounded first derivatives.

(ii) If the system considered is   a gradient  system, then

\begin{equation}
\delta P_t\psi=\heatsemi\psi,
\label{eq: gradient 1}
\end{equation}
for all  bounded  $C^2$ q-forms $\psi$, 
 provided that  $E\left(\sup_{s\le t} |TF_s|^q\right)$ is finite for each
 $t>0$. 

 In particular
these hold if $|\nabla X|$ is bounded and $H_1$ is bounded above.

\end{proposition*}

However the following often has advantage when $P_tf$ is known to be $BC^1$.

\begin{proposition*}
 Assume $\A=\half \triangle+L_Z$. 
Let $M$ be a complete Riemannian manifold with Ricci curvature bounded
from below by $-c(1+\rho^2(x))$. 
Suppose $d\rho(Z(x))\le c[1+\rho(x)]$ and $H_{1+\delta}(x)(v,v)
\le c\ln[1+\rho(x)]|v|^2$ for all $x$ and $v$. Here
$c$,  and $\delta>0$ are constants. Then
$$dP_tf=\delta P_t(df)$$
for all $f$ in $C_K^\infty$ provided $d(P_tf)$
is bounded uniformly in each $[0,T]$. 
\end{proposition*}

\noindent
{\bf Proof:}
Let $\phi_0$ be a bounded $C^2$ 1-form. We shall show that a solution 
$P_t\phi$ to ${\partial\phi_t\over \partial t}=\triangle^1\phi_t+L_Z\phi_t$ 
starting from $\phi_0$ and bounded on $[0,T]\times M$ is given by 
$E\phi_0(v_t)$ and then note
 that $d(P_tf)=P_t(df)$ for smooth functions to finish the proof.
Let $\tau_n(x_0)$ be the first exit time of $F_t(x_0)$ from the ball $B(n)$
radius $n$, centred at $p$. Since $P_t\phi$ is smooth, we apply 
It\^o's formula to get:
$$P_{T-t}\phi(v_t)=P_T\phi(v_0)+\int_0^t\nabla P_{T-s}\phi( XdB_s)
+\int_0^t P_{T-s}\phi(\nabla X(v_s)dB_s).$$
Replace $t$ by $t\wedge \tau^n$ in the above inequality to obtain:
\begin{eqnarray*}
P_{T-t\wedge \tau^n}\phi(v_{t\wedge \tau^n})&=&P_T\phi(v_0)
+\int_0^{t\wedge \tau^n}\nabla P_{T-s}\phi( XdB_s)\\
&&+\int_0^{t\wedge \tau^n} P_{T-s}\phi(\nabla X(v_s)dB_s).  \end{eqnarray*}
This gives:
\begin{equation}\label{eq: appendix 30}
E\phi(v_T)\chi_{T\le \tau^n} +EP_{T-\tau^n}\phi(v_{\tau^n})\chi_{\tau^n<T}
=P_T\phi(v_0).\end{equation}
But under the condition $H_{1+\delta}(x)\le c\ln[1+\rho(x)]$,
\begin{equation}\label{eq: appendix 10}
E|v_{\tau^n}|^{1+\delta}\chi_{\tau^n<T}\le \exp^{C_nT/2}.\end{equation}
Here $C_n=\sup_{x\in B(n)}\sup_{|v|\le 1}H_{1+\delta}(x)(v,v)\le c\ln(1+n)$. 
See  \cite{flow} for details.  On the other hand \cite{flow}, there is a 
constant $k_0>0$  such that for each  $\beta>0$

\begin{equation}
P(\tau^n(x)<T)\le {1\over n^\beta}[1+\rho(x)]^\beta \exp^{k_0[1+\beta^2]T}.
\label{eq: appendix 20}\end{equation}
Take numbers  $\delta'>0$, and  $p>1$, $q>1$ such that
${1\over p}+{1\over q}=1$ and $p(1+\delta')=1+\delta$. Then
$$\sup_n E|P_{T-\tau^n} \phi(v_{\tau^n})\chi_{\tau^n<T}|^{1+\delta'}\le
k\sup_n \left[E|v_{\tau^n}\chi_{\tau^n<T}|^{p(1+\delta')}\right]^{1\over p}
\left[P(\tau^n<T)\right]^{1\over q}.$$
Here $k$ is a constant. We have used the assumption that $P_t\phi$ is 
uniformly bounded on $[0,T]$.

  By choosing $\beta$ sufficiently big, we see, from (\ref{eq: appendix 10})
 and (\ref{eq: appendix 20}), that the right hand side of the inequality  is
 finite. Thus
$|P_{T-\tau^n}\phi(v_{\tau^n})\chi_{\tau^n<T}|$ is uniformly integrable.
Passing to the limit $n\to \infty$ in (\ref{eq: appendix 30}), we have shown:
$$E\phi(v_T)=P_T\phi(v_0).$$
\hfill\rule{3mm}{3mm}
\bigskip

There are also parallel results for higher order forms.

\begin{proposition*}\label{pr: 7.8}
 Suppose the s.d.e. (\ref{eq: 0}) is strongly 1-complete and
 $T_xF_t$ is also strongly 1-complete. Let  $f\in BC^2$, then
\begin{equation}\label{eq: second derivative}
\begin{array}{lll}\nabla d(P_tf)(u,v) &=& E\nabla (df)(T_xF_t(u), T_xF_t(v))
\\ &+&Edf(\nabla (TF_t)(u,v))\end{array}\end{equation}
for all $u, v\in T_xM$, if for each $t>0$ and compact set $K$, there is a 
constant $\delta>0$ such that

\begin{equation}\label{second derivative 1}
\sup_{x\in K} E|T_xF_t|^{2+\delta}<\infty,\end{equation}
 and

\begin{equation}\label{second derivative 2}
\sup_{x\in K} E|\nabla T_xF_t|^{1+\delta}<\infty.\end{equation}
In particular (\ref{eq: second derivative}) holds if the first three
 derivatives  of $X$ and the first two derivatives of $A$ are bounded.

\end{proposition*}

\noindent{\bf Proof:}
First $dP_tf=\delta P_t(df)$ from a result in \cite{flow}.
Let $u,v\in T_xM$.  Take a  smooth map $\sigma_1\colon [0,s_0] \to M$ such that
$\dot\sigma_1(0)=u$. Let $v(s)\in T_{\sigma_1(s)}M$ be the parallel
translate of $v$ along $\sigma_1$. Suppose its image is contained in a compact
 set $K$. Then 
$df_{F_t(\sigma_1(s))}\left(T_{\sigma_1(s)} F_t(v(s))\right)$
 is a.s. differentiable in $s$ for each $t>0$.  So for almost all $\omega$,

\begin{eqnarray*}
I_s
&=&{df_{F_t(\sigma_1(s))}\left(T_{\sigma_1(s)} F_t(v(s))\right)
 -df\left(T_xF_{t}(v)\right)\over s}\\
&=&{1\over s}\int_0^s {\partial\over \partial r}\left[df
\left(T_{(\sigma_1(r))}F_t(v(r))\right)\right]dr\\
&=&{1\over s}\int_0^s \nabla df\left(TF_t(\dot \sigma_1(r)),
 TF_t(v(r))\right)dr+
{1\over s}\int_0^s df\left(\nabla TF_t(\dot \sigma_1(r), v(r))\right)dr\\
\end{eqnarray*}

But the integrand  of the right hand side is continuous in $r$ in $L_1$, so
$E\lim_{s\to 0}I_s=\lim_{s\to 0}EI_s$. Thus
\begin{eqnarray*}&& \nabla d(P_tf)(u,v)=\nabla (\delta P_t(df))(u,v)\\
&=&E\nabla (df)(TF_t(u), TF_t(v))
 +Edf(\nabla (TF_t)(u,v)).\end{eqnarray*}
For the last part observe that if the s.d.e. is strongly 2-complete, then 
$TF_t$ is strongly 1-complete and apply lemma A\ref{le: 7.2}. 
\hfill\rule{3mm}{3mm}

\bigskip

For elliptic systems, we may use the previous weak derivatives argument.
Just notice that for two $C_K^\infty$ vector fields $\Lambda_1$ and
$\Lambda_2$,
$$L_{\Lambda_2}\L_{\Lambda_1}(P_tf)(x)=
\nabla^2 P_tf(x)(\Lambda_2(x), \Lambda_1(x))+
<\nabla P_tf(x),\nabla \Lambda_1(\Lambda_2(x))>_x$$
and 
\begin{eqnarray*}
\L_{\Lambda_2}df\circ T_xF_t(\Lambda_1(x))&=&
\nabla df(T_xF_t(\Lambda_2), T_xF_t(\Lambda_1)\\
&&+ df\circ \nabla T_xF_t(\Lambda_2, \Lambda_1)+
df\circ T_xF_t(\nabla \Lambda_1(\Lambda_2(x))).
\end{eqnarray*}
In this case  the number $\delta$ in the assumption can be taken to be zero,
  but the required equality (\ref{eq: second derivative}) holds only
 almost surely.  However this  is usually  enough for our purposes.


\begin{thebibliography}{10}

\bibitem{Bismut}
J.~Bismut.
\newblock {\em Large deviations and the {M}alliavin calculus}.
\newblock Birkhauser, 1984.

\bibitem{CH73}
P.~R. Chernoff.
\newblock Essential self-adjointness of powers of generators of hyperbolic
  equations.
\newblock {\em Journal of Functional Analysis}, 12:401--414, (1973).

\bibitem{DA-EL-ZA}
G.~DaPrato, D.~Elworthy, and J.~Zabczyk.
\newblock Strong {F}eller property for stochastic semilinear equations.
\newblock Warwick preprints, 1/1992. To appear in {S}tochastic {A}nalysis and
  its Applications.

\bibitem{ell-koh}
R.~Elliott and M.~Kohlmann.
\newblock Integration by parts, homogeneous chaos expansions and smooth
  densities.
\newblock {\em Annals of Probability}, 17(1):194--207, 1989.

\bibitem{ELflour}
K.~D. Elworthy.
\newblock Geometric aspects of diffusions on manifolds.
\newblock In P.~L. Hennequin, editor, {\em Ecole d'Et\'e de Probabilit\'es de
  Saint-Flour XV-XVII, 1985, 1987. Lecture Notes in Mathematics}, volume 1362,
  pages 276--425. Springer-Verlag, 1988.

\bibitem{ELflow}
K.~D. Elworthy.
\newblock Stochastic flows on {R}iemannian manifolds.
\newblock In M.~A. Pinsky and V.~Wihstutz, editors, {\em Diffusion processes
  and related problems in analysis, volume II. {B}irkhauser Progress in
  Probability}, pages 37--72. Birkhauser, Boston, 1992.

\bibitem{EL-YOR}
K.~D. Elworthy and M.~Yor.
\newblock Conditional expectations for derivatives of certain stochastic flows.
\newblock In J.~Az\'ema, P.A. Meyer, and M.~Yor, editors, {\em Sem. de Prob.
  XXVII. Lecture Notes in Maths. 1557}, pages 159--172. Springer-Verlag, 1993.

\bibitem{ELbook}
K.D. Elworthy.
\newblock {\em Stochastic {D}ifferential {E}quations on {M}anifolds}.
\newblock Lecture Notes Series 70, Cambridge University Press, 1982.

\bibitem{Freidlin85}
M.~Freidlin.
\newblock {\em Functional {I}ntegration and {P}artial {D}ifferential
  {E}quations}.
\newblock Princeton Unviersity Press, 1985.

\bibitem{Krylov93}
N.V. Krylov.
\newblock Quasiderivatives for solutions of {I}t\^o's stochastic equations and
  their applications.
\newblock Preprint.

\bibitem{application}
Xue-Mei  Li.
\newblock Stochastic differential equations on noncompact manifolds: moment
  stability and its topological consequences.
\newblock Warwick Preprints: 15/1994. To appear in Probab. Theory Relat.
  Fields.

\bibitem{flow}
Xue-Mei  Li.
\newblock Strong p-completeness and the existence of smooth flows on noncompact
  manifolds.
\newblock Warwick preprints: 28/1993. To appear in Probab. Theory Relat.
  Fields.

\bibitem{Li.thesis}
Xue-Mei  Li.
\newblock Stochastic {F}lows on {N}oncompact {M}anifolds.
\newblock Ph.D. thesis, University of Warwick., 1992.

\bibitem{Li-Zhao}
Xue-Mei Li and H.Z. Zhao.
\newblock Gradient estimates and the smooth convergence of approximate
              travelling waves for reaction-diffusion equations. In preparation.

\bibitem{Norris93}
J.~Norris.
\newblock Path integral formulae for heat kernels and their derivatives.
\newblock {\em Probability Theory and Related Fields}, 94:525--541, 1993.

\bibitem{PE-ZA}
S.~Peszat and J.~Zabczyk.
\newblock Strong {F}eller property and irreducibility for diffusions on
  {H}ilbert spaces.
\newblock Preprint of {P}olish {A}cademy of {S}ciences, 1993.

\bibitem{Schwartz66}
L.~Schwartz.
\newblock {\em Th\'eorie des distributions}.
\newblock Hermann, 1966.

\bibitem{Watanabe}
S.~Watanabe.
\newblock {\em Lectures on stochastic differential equations and {M}alliavin
  Calculus.}
\newblock Bombay: {T}ata {I}nstitute of {F}undamental {R}esearch,
  Springer-Verlag, 1984.

\end{thebibliography}

Address:\\
Mathematics Institute, \\University of Warwick,\\ Coventry CV4, 7AL, U.K.
\end{document}